# A scale-dependent finite difference method for time fractional derivative relaxation type equations

XiaoTing Liu[a], HongGuang Sun[a,*], Yong Zhang[a,b], Zhuojia Fu[a]

a Institute of Soft Matter Mechanics, College of Mechanics and Materials, Hohai University, Nanjing, Jiangsu, 210098, China

b Department of Geological Sciences, University of Alabama, Tuscaloosa, AL 35487,USA

**Abstract**

Fractional derivative relaxation type equations (FREs) including fractional diffusion equation and fractional relaxation equation, have been widely used to describe anomalous phenomena in physics. To utilize the characteristics of fractional dynamic systems, this paper proposes a scale-dependent finite difference method (S-FDM) in which the non-uniform mesh depends on the time fractional derivative order of FRE. The purpose is to establish a stable numerical method with low computation cost for FREs by making a bridge between the fractional derivative order and space-time discretization steps. The proposed method is proved to be unconditional stable with $(2-\alpha)$-th convergence rate. Moreover, three examples are carried out to make a comparison among the uniform difference method, common non-uniform method and S-FDM in term of accuracy, convergence rate and computational costs. It has been confirmed that the S-FDM method owns obvious advantages in computational efficiency compared with uniform mesh method, especially for long-time range computation (e.g. the CPU time of S-FDM is ~1/400 of uniform mesh method with better relative error for time T=500 and fractional derivative order $\alpha$=0.4)..

1. Introduction

Fractional derivative which is almost as old as integer-order derivative[1, 2], has a comprehensive applications in hydrology[3, 4, 31-33], image processing[5-7], finance[8], physics[9], etc. Previous investigations have confirmed that, fractional-order derivative models are more suitable to establish and develop dynamic system on account of its history dependence or/and non-local/long-range interaction. According to the discrepancy of natural physical processes and fractional definitions, there are different kinds of fractional derivative models to describe various anomalous phenomena and dynamic systems[10-13, 34, 37].

In the past decades, there are many theoretical and numerical algorithm for FREs have drawn attention and attracted increasing interests. Since there are so many papers devoted to this field, here we just review several previous literatures[14-25, 31-52]. Podlubny[24] proved the fractional differential equations to model anomalous diffusion phenomena. Meerschaert et al.[25] developed a finite difference approximations for fractional flow equations with variable coefficients on a finite domain. K.Diethelm et al.[27] made a great effort to provide the predict-corrected method for fractional Relaxation-Oscillation equation. But we should denote that though numerical methods with same discretization step is easily gained and calculated, the computation efficiency based on uniform step is a remarkable problem in applications[31-51]. Therefore, the high-efficient numerical methods for FREs are still necessary for the physical modeling and experiment results analysis in consideration of high computational cost for long-range computation owing to its long-time memory and spatial non-locality. Hence, non-uniform mesh method as a promising approach, has been proposed with ambition to improve computation efficiency, due to its distinct advantages for accurately solving FREs with less mesh points[28,29]. To our knowledge, Li et al. [28] presented a finite difference method to deal with nonlinear fractional differential equations by using non-uniform meshes with non-equidistant step-size. Zhang et al. [14] extended the finite difference method for time fractional diffusion equation equation on non-uniform meshes. However, one can notice that the non-uniform collocation methods in previous literatures have no relationship with fractional order of FREs, while the behaviors of fractional relaxation systems have close relationship with their differential orders. Thus, how to determine the non-uniform mesh based on the characteristics of fractional relaxation systems is an interesting and important topic to improve the computational efficiency of numerical method.

The purpose of this paper is to introduce a non-uniform mesh method based on the relationship of fractional derivative order and strategy of non-uniform collocation. Chen et al.[16, 17, 33] introduced a new approach of space-time metrics and further verified it through a series of applications. Here we borrow the idea of Hausdroff metric to establish the relationship between strategy of non-uniform mesh and fractional derivative order. For illustration purpose, we only consider the time derivative FREs which account for the heavy tail phenomena caused by historical memory, to demonstrate its efficiency and accuracy. Since the relationship between intensity of the memory property and time fractional order could be called as scale-dependence phenomenon from application approach, we name this method as scale-dependent finite difference method (S-FDM). To make our illustration clear and brief in the sequel, we introduce two concepts named clock time in Euclidean metric and power-law time in Hausdroff metric. After that, we could illustrate that the proposed method could be regard as the uniform mesh on the power-law time while it shows non-uniform mesh property on the clock time. Meanwhile, what is important is that it combines the two advantages of non-uniform mesh methods on high accuracy and uniform mesh methods on low computation cost.

The rest of this paper is organized as follows. In section 2, the stability and convergence of numerical scheme are proved. Then numerical results on different schemes are calculated to show the property of scale-dependent finite difference method in section 3. Then, further comparison and discussion on presented method are provided to verify our statement in section 4. Finally, some conclusions are drawn in section 5.

2. **Theory and Methodologies**

Consider the following time FRE (could be relaxation/diffusion equation)

$$^C_0D^\alpha_t u(x,t) = Ku(x,t) + f(x,t) \tag{1}$$

with the initial and boundary conditions

$$\begin{aligned} u(x,t)|_{\partial\Omega} &= 0, \quad 0 < t \leq T, \\ u(x,0) &= w(x,t), \quad x \in \partial\Omega, \end{aligned} \tag{2}$$

where the $f$ is the source term and K represents relaxation parameter. The well-known Caputo fractional operator is

$$^C_0D^\alpha_t u(x,t) = \int_0^t f'(\xi)(t-\xi)^{-\alpha} d\xi, \quad 0 < \alpha < 1$$

Here, we separate the interval into integer N subintervals with $0 = t_0 < t_1 < \cdots < t_N = T$. The time steps are denoted as $\tau_k = t_k - t_{k-1}$.

We divide equidistant step-sizes $\Delta t_\alpha = T^\alpha/N$ from the power-law time interval $[0, T^\alpha]$ with fractional order α and integer N. If the order α is equal to 1, the power-law time would become clock time due to the relationship between the Hausdroff metric and Euclidean metric. For the sake of intuitionistic exhibition, we draw the Figure.4 to show the relationship of the collocation modes between Hausdroff metric and Euclidean metric. It is unambiguous to show the relative positions between the points on power-law time (Hausdroff metric γ=0.5) with the equidistant step and the parallel positions on clock time when collocation nodes are 50. Hence, the presented mesh is uniform on power-law time scale as well as non-uniform in the clock time scale.

Then we have $^n_1\tau_{max} = \max_{1\leq k \leq n} \tau_k = \tau_n$ and $^n_1\tau_{min} = \min_{1\leq k \leq n} \tau_k = \tau_1$ when α ∈ (0, 1].

For the proposed meshes above, it has the following result [14].

**Lemma 2.1.** For $0 < \alpha < 1$ and $f(t) \in C^2, t \in [0,T]$, n is any integer in [0, N]. Let $\Delta t_\alpha = T^\alpha / N$. It holds that

$$\int_0^{t_n} f'(\xi)(t_n - \xi)^{-\alpha} d\xi = \sum_{k=1}^n \frac{f(t_k) - f(t_{k-1})}{\tau_k} \int_{t_{k-1}}^{t_k} (t_n - \xi)^{-\alpha} d\xi + R^n \tag{3}$$

$$|R^n| \leq \frac{\tau_n^{2-\alpha}}{2(1-\alpha)} \max_{t_{n-1} \leq t \leq t_n} |f''(t)| + \frac{\alpha}{8} \sum_{k=1}^{n-1} \max_{t_{k-1} \leq t \leq t_k} |f''(t)| \tau_k^2 \int_{t_{k-1}}^{t_k} (t_n - \xi)^{-\alpha-1} d\xi$$

$$\leq (\frac{1}{2(1-\alpha)} + \frac{1}{8}) \tau_n^{2-\alpha} \max_{0 \leq t \leq t_n} |f''(t)| \tag{4}$$

**Proof.** The integral could be written as

$$\int_0^{t_n} f'(\xi)(t_n - \xi)^{-\alpha} d\xi = \int_0^{t_{n-1}} f'(\xi)(t_n - \xi)^{-\alpha} d\xi + \int_{t_{n-1}}^{t_n} f'(\xi)(t_n - \xi)^{-\alpha} d\xi$$

and deal with the formula (4) by two parts. Firstly, utilizing linear interpolation of $f$, we obtain

$$\int_0^{t_{n-1}} f'(\xi)(t_n - \xi)^{-\alpha} d\xi$$

$$= \left[ (t_n - \xi)^{-\alpha} f(\xi) \right]_0^{t_{n-1}} - \alpha \int_0^{t_{n-1}} f(\xi)(t_n - \xi)^{-\alpha-1} d\xi$$

$$= \tau_n^{-\alpha} f(t_{n-1}) - t_n^{-\alpha} f(0) - \alpha \sum_{k=1}^{n-1} \int_{t_{k-1}}^{t_k} \frac{(t_k - \xi) f(t_{k-1}) + (\xi - t_{k-1}) f(t_k)}{\tau_k} (t_n - \xi)^{-\alpha-1} d\xi - (R_1)^n$$

$$= \tau_n^{-\alpha} f(t_{n-1}) - t_n^{-\alpha} f(0) + \sum_{k=1}^{n-1} f(t_{k-1}) (t_n - t_{k-1})^{-\alpha} - \sum_{k=1}^{n-1} f(t_k)(t_n - t_k)^{-\alpha}$$

$$+ \sum_{k=1}^{n-1} \frac{f(t_k) - f(t_{k-1})}{\tau_k} \int_{t_{k-1}}^{t_k} (t_n - \xi)^{-\alpha} d\xi - (R_1)^n$$

$$= \sum_{k=1}^{n-1} \frac{f(t_k) - f(t_{k-1})}{\tau_k} \int_{t_{k-1}}^{t_k} (t_n - \xi)^{-\alpha} d\xi - (R_1)^n$$

where

$$(R_1)^n = \alpha \sum_{k=1}^{n-1} \int_{t_{k-1}}^{t_k} \frac{1}{2} f''(\xi_k)(\xi - t_k)(\xi - t_{k-1})(t_n - \xi)^{-\alpha-1} d\xi, \quad \xi_k \in (t_{k-1}, t_k)$$

$$\begin{aligned}
|(R_1)^n| &\leq \frac{\alpha}{8} \sum_{k=1}^{n-1} \tau_k^2 |f''(t_k)| \int_{t_{k-1}}^{t_k} (t_n - \xi)^{-\alpha-1} d\xi \\
&= \frac{1}{8} \sum_{k=1}^{n-1} \tau_k^2 |f''(t_k)| ((t_n - t_{k-1})^{-\alpha} - (t_n - t_k)^{-\alpha}) \\
&\leq \frac{\alpha \, _1^{(n-1)} \tau_{\max}^2}{8} \max_{0 \leq t \leq t_{n-1}} |f''(t_k)| \int_0^{t_{k-1}} (t_n - \xi)^{-\alpha-1} d\xi \\
&= \frac{_1^{(n-1)} \tau_{\max}^2}{8} \max_{0 \leq t \leq t_{n-1}} |f''(t_k)| (\tau_n^{-\alpha} - t_n^{-\alpha})
\end{aligned} \quad (5)$$

Considering equidistant step-sizes $\Delta t_\alpha = T^\alpha / N$ on the interval $[0, T^\alpha]$, we get

$$\begin{aligned}
|(R_1)^n| &\leq \frac{\alpha \tau_{n-1}^2}{8} \max_{0 \leq t \leq t_{n-1}} |f''(t_k)| \int_0^{t_{k-1}} (t_n - \xi)^{-\alpha-1} d\xi \\
&= \frac{\tau_{n-1}^2}{8} \max_{0 \leq t \leq t_{n-1}} |f''(t_k)| (\tau_n^{-\alpha} - t_n^{-\alpha}) \\
&\leq \frac{1}{8} \max_{0 \leq t \leq t_{n-1}} |f''(t_k)| ((((k-1)^{1/\alpha} - k^{1/\alpha}) \Delta t_\alpha^{1/\alpha})^{2-\alpha} - t_n^{-\alpha} (((k-2)^{1/\alpha} - (k-1)^{1/\alpha}) \Delta t_\alpha^{1/\alpha})^2) \\
&\leq \frac{1}{8} \max_{0 \leq t \leq t_{n-1}} |f''(t_k)| (\tau_n^{2-\alpha} - t_n^{-\alpha} \tau_n^2) \\
&\leq \frac{1}{8} \max_{0 \leq t \leq t_{n-1}} |f''(t_k)| \tau_n^{2-\alpha}
\end{aligned} \quad (6)$$

Furthermore, with Taylor expansion, the error in interval $[t_{n-1}, t_n]$ could be

$$\begin{aligned}
|(R_2)^n| &= \left| \int_{t_{n-1}}^{t_n} f'(\xi)(t_n - \xi)^{-\alpha} d\xi - \int_{t_{n-1}}^{t_n} \frac{f(t_n) - f(t_{n-1})}{\tau_n} (t_n - \xi)^{-\alpha} d\xi \right| \\
&\leq \frac{\tau_n^2}{2} \max_{t_{n-1} \leq t \leq t_n} |f''(t)| \frac{\tau_n^{-\alpha}}{1-\alpha} \\
&\leq \frac{\tau_n^{2-\alpha}}{2(1-\alpha)} \max_{t_{n-1} \leq t \leq t_n} |f''(t)|
\end{aligned} \quad (7)$$

Combining (6) and (7), we complete the proof.

### 2.1 Derivation of the difference scheme

Let $\Delta_t u(x, t_k) = u(x, t_{k+1}) - u(x, t_k)$, and $db_{k,j} = k^{1/\alpha} - j^{1/\alpha}$

The scale-dependent finite difference numerical scheme for Caputo time fractional derivative is written as follows:

$$\frac{\partial^\alpha u(x,t_{n+1})}{\partial t^\alpha} \approx \frac{1}{\Gamma(1-\alpha)}\sum_{k=1}^{n}\frac{\Delta_t u(x,t_k)}{\tau_{k+1}}\int_{t_k}^{t_{k+1}}\frac{d\xi}{(t_{n+1}-\xi)^\alpha}$$

$$=\frac{1}{\Gamma(2-\alpha)}\sum_{k=1}^{n}\frac{\Delta_t u(x,t_k)}{\tau_{k+1}}[(t_{n+1}^{1/\alpha}-t_k^{1/\alpha})^{1-\alpha}-(t_{n+1}^{1/\alpha}-t_{k+1}^{1/\alpha})^{1-\alpha}] \qquad (8)$$

$$=\frac{-\Delta t_\alpha^{-1}}{\Gamma(2-\alpha)}\sum_{k=1}^{n}\frac{\Delta_t u(x,t_k)}{b_{k+1,k}}[(b_{n+1,k+1})^{1-\alpha}-(b_{n+1,k})^{1-\alpha}]$$

where $\tau_k = t_k - t_{k-1}$. The $u(x,t)$ is the solution of the equation(1).

Here, we denote

$$\chi_k^n = \frac{1}{\tau_{k+1}}\int_{t_k}^{t_{k+1}}\frac{d\xi}{(t_{n+1}-\xi)^\alpha}$$
$$= \frac{1}{(1-\alpha)\tau_{k+1}}[(t_{n+1}^{1/\alpha}-t_k^{1/\alpha})^{1-\alpha}-(t_{n+1}^{1/\alpha}-t_{k+1}^{1/\alpha})^{1-\alpha}], \qquad 0 \le k \le n \qquad (9)$$

Then we set $\chi_k^n = \frac{1}{\Gamma(1-\alpha)}(t_{n+1}-\xi_k)^{-\alpha}$, $\xi_k \in [t_{k-1},t_k]$,

Due to the monotonic increasing function $(t-s)^{-\alpha}$, it is easy to get the inequality

$$\chi_n^n \ge \chi_{n-1}^n \ge \cdots \ge \chi_k^n \ge \cdots \ge \chi_0^n \ge 0 \qquad (10)$$

while, it is also noticed that

$$\chi_0^n > t_{n+1}^{-\alpha} > T^{-\alpha} \qquad (11)$$

Consider Eqs.(1) and (2) in interval [0,T] and Eq.(8), the numerical function could be

$$\frac{1}{\Gamma(1-\alpha)}\sum_{k=0}^{n}\Delta_t u(x,t_k)\chi_k^n = Ku(x,t_{n+1}) + f(x,t_{n+1}) + R(t_{n+1}),\ 0 \le t_{n+1} \le T$$
$$u(x,t_{n+1})|_{\partial\Omega}=0, \quad 0 < t_{n+1} \le T, \qquad (12)$$
$$u(x,0) = w(x), \quad x \in \partial\Omega,$$

where $R(t_{n+1})$ is the truncation errors of moment $t_{n+1}$. Obviously, we could omit the small error $R(t_{n+1})$ and the $\bar{u}(x,t_{n+1})$ is numerical approximation of the exact solution. After that, we have following discrete difference scheme

$$\frac{1}{\Gamma(1-\alpha)}\sum_{k=0}^{n}\Delta_t \bar{u}(x,t_k)\chi_k^n = K\bar{u}(\mathbf{x},t_{n+1}) + f(x,t_{n+1}),\ 0 \le t_{n+1} \le T \qquad (13)$$

with initial and boundary condition

$\bar{u}(x, t_{n+1})|_{\partial\Omega} = 0, \quad 0 < t_{n+1} \leq T,$

$\bar{u}(x, 0) = w(x), \quad x \in \partial\Omega,$

**Stability and convergence**

Firstly, the formulas below are introduced

$\langle m, n \rangle = \int_\Omega mn d\mathbf{x}$

$\|m\| = \sqrt{\langle m, m \rangle}$

$\|\nabla m\| = \sqrt{\langle \nabla m, \nabla m \rangle}$

$\|m\|_1 = \sqrt{\langle \|m\|^2, \|\nabla m\|^2 \rangle}$

Then the Green's first formula is

$$\langle m, \Delta w \rangle + \langle \nabla m, \nabla w \rangle = \oint_{\partial\Omega} m(\nabla m \cdot \vec{n}) dS \tag{14}$$

where $\Delta$ and $\nabla$ represent Laplace and gradient operator separately. $\vec{n}$ is the outer normal vector of limited surface dS.

**Lemma 2.2.** For the Eqs.(1) and (2), the stability of the difference scheme (12) is unconditional

$$\|\nabla \bar{u}(x, t_{n+1})\|^2 \leq \|\nabla w(x)\|^2 + \frac{T^\alpha \Gamma(1-\alpha)}{2} \max_{1 \leq k \leq n+1} \|f(x, t_k)\|^2, \quad 0 \leq t_n \leq T$$

**Proof.** The above Eq.(13) takes the inner product with $-2\nabla \bar{u}(x, t_{n+1})$. Utilizing formula (14) and denoting the zero boundary, we get

$\frac{2}{\Gamma(1-\alpha)} \chi_n^n \|\nabla \bar{u}^{n+1}\|^2 + 2\|\Delta \bar{u}^{n+1}\|$

$= \frac{1}{\Gamma(1-\alpha)} \left( 2 \sum_{k=0}^{n} (\chi_{k+1}^n - \chi_k^n) \langle \nabla \bar{u}^k, \nabla \bar{u}^{n+1} \rangle + 2\chi_0^n \langle \nabla \bar{u}^0, \nabla \bar{u}^{n+1} \rangle \right) - 2\langle f^{n+1}, \Delta \bar{u}^{n+1} \rangle$

$\leq \frac{1}{\Gamma(1-\alpha)} \left( \sum_{k=0}^{n} (\chi_{k+1}^n - \chi_k^n)(\|\nabla \bar{u}^k\|^2 + \|\nabla \bar{u}^{n+1}\|^2) + \chi_0^n (\|\nabla \bar{u}^0\|^2 + \|\nabla \bar{u}^{n+1}\|^2) \right) + \frac{1}{2}\|f^{n+1}\|^2 + 2\|\Delta \bar{u}^{n+1}\|^2$

$= \frac{1}{\Gamma(1-\alpha)} \left( \sum_{k=0}^{n} (\chi_{k+1}^n - \chi_k^n)\|\nabla \bar{u}^k\|^2 + \chi_n^n \|\nabla \bar{u}^n\|^2 + \chi_0^n \|\nabla \bar{u}^0\|^2 \right) + \frac{1}{2}\|f^{n+1}\|^2 + 2\|\Delta \bar{u}^{n+1}\|^2$

Considering Eqs. (10) and (11), it could be gain

$$\chi_n^n \|\nabla \bar{u}^{n+1}\|^2 \leq \sum_{k=0}^{n} (\chi_{k+1}^n - \chi_k^n)\|\nabla \bar{u}^k\|^2 + \chi_0^n \left( \|\nabla \bar{u}^0\|^2 + \frac{T^\alpha \Gamma(1-\alpha)}{2}\|f^{n+1}\|^2 \right), 0 \leq t_{n+1} \leq T$$

The proof have been completed    □

using mathematic induction, we have

$$\left\|\nabla \bar{u}^{n+1}\right\|^2 \leq \left\|\nabla \bar{u}^{0}\right\|^2 + \frac{T^\alpha \Gamma(1-\alpha)}{2}\left\|f^{n+1}\right\|^2, 0 \leq t_{n+1} \leq T$$

then

$$\chi_n^n \left\|\nabla \bar{u}^{n+1}\right\|^2 \leq \sum_{k=0}^{n}(\chi_{k+1}^n - \chi_k^n)\left\|\nabla \bar{u}^{-k}\right\|^2 + \chi_0^n \left(\left\|\nabla \bar{u}^{0}\right\|^2 + \frac{T^\alpha \Gamma(1-\alpha)}{2}\left\|f^{n+1}\right\|^2\right)$$

$$\leq \sum_{k=0}^{n}(\chi_{k+1}^n - \chi_k^n)\left(\left\|\nabla \bar{u}^{0}\right\|^2 + \frac{T^\alpha \Gamma(1-\alpha)}{2}\left\|f^{n+1}\right\|^2\right) + \chi_0^n \left(\left\|\nabla \bar{u}^{0}\right\|^2 + \frac{T^\alpha \Gamma(1-\alpha)}{2}\left\|f^{n+1}\right\|^2\right)$$

$$= \chi_n^n \left(\left\|\nabla \bar{u}^{0}\right\|^2 + \frac{T^\alpha \Gamma(1-\alpha)}{2}\left\|f^{n+1}\right\|^2\right), \qquad 0 \leq t_{n+1} \leq T$$

Set

$$e(x, t_{n+1}) = u(x, t_{n+1}) - \bar{u}(x, t_{n+1}),\ 0 \leq t_{n+1} \leq T$$

Noting that the error of eq.(12), we have the below equation

$$\frac{1}{\Gamma(1-\alpha)}\sum_{k=0}^{n}\Delta_t e(x,t_k)\chi_k^n = Ke(x,t_{n+1}) + R(t_{n+1}),\ 0 \leq t_{n+1} \leq T$$
$$e(x,t_{n+1})|_{\partial\Omega} = 0,\ \ 0 < t_{n+1} \leq T, \qquad (15)$$
$$e(x,0) = 0,\ \ x \in \partial\Omega,$$

The truncation error $R(t_{n+1})$ is close related to discrete Caputo derivative part and diffusion part. Thus, with

Lemma 2.2, we gain

$$\left\|\nabla e(x,t_{n+1})\right\|^2 \leq \frac{T^\alpha \Gamma(1-\alpha)}{2}\left\|R(t_{n+1})\right\|^2 \leq \frac{T^\alpha \Gamma(1-\alpha)}{2}\max_{0 \leq t_k \leq T}\left\|R(t_k)\right\|^2$$
$$\leq C\tau_n^{2-\alpha} \qquad 0 \leq t_{n+1} \leq T$$

**Theorem 2.3** Assuming that $u(x,t) \in \mathcal{C}_{x,t}^{2,2}$, we have

$$\left\|\nabla u(x,t_{n+1}) - \nabla \bar{u}(x,t_{n+1})\right\| \leq \sqrt{\frac{T^\alpha \Gamma(1-\alpha)}{2}}C\tau_n^{2-\alpha}, \quad 0 \leq t_{n+1} \leq T$$

Applying the *Poincaré inequality* with **Theorem 2.1,** the numerical solution of fractional diffusion equation(13) is convergence in $H^1$ norm as $\tau_n \to 0$.

## 3. Numerical Experiments

In this section, we carry out numerical experiments for a variety of time FREs by the proposed numerical difference method call S-FDM. All tests implement on the condition of MATLAB R2012a with a desktop computer (Lenovo yangtian T4990v=00) having the following configuration: Intel(R) Core(TM) i5-4590 CPU 3.30 GHz and 16.00G RAM. In the paper, CPU time is the computation cost on mentioned computer condition.

**Example1.** Consider the following fractional relaxation equation

$$\begin{cases} {}_0^C D_t^\gamma u(t) + Bu(t) = 0 \\ u(0) = A \end{cases}, \quad t \in [0,T] \tag{16}$$

where we gain the order $\gamma = 0.5$, initial value A=10 on the numerical experiments. The exact solution is $uexact(t) = AE_\gamma(-Bt^\gamma)$.

On Figure 1, the conditions of numerical computation are adopted with same collocations 100 and equivalent T=10. The curves show the results of different relaxation coefficient B=3, 4, 5 with rectangular formula on clock time with $\Delta t = T/N$ and power-law time with $\Delta t_\alpha = T^\alpha/N$, separately. The curves on clock time show the numerical results are not stable apparently, especially in the initial time. Here, the numerical results for FRE(16) become unstable significantly while the relaxation coefficient B is equal to 5. In contrast, the curves with S-FDM behave better that could obtain convergent numerical solutions with total mentioned coefficients B. Therefore, the results reflect that the rectangular formula with S-FDM could fit the wide relaxation rate better at the same computation condition.

When the γ=0.5, B=1 and A=10, the exact solution of FRE(16) could be $uexact(t) = 10e^t erfc(\sqrt{t})$. From the Table1, the implicit numerical methods with the two methods are compared in terms of accuracy, rate of convergence and computation cost. It is obviously shown that the S-FDM just obtains same accuracy with one quarter collocations while the implicit uniform FDM needs 400 collocations. What is more, the scale-dependent difference method has a super-linear convergence rate that the rate on Table 1 is similar to 1.5 that is equal to (2-γ) while the uniform FDM get first order convergence rate. It means that the S-FDM reduces the computation costs 4 times on the Hausdroff metric $\gamma = 0.5$. Therefore, it is clear that the S-FDM on Hausdroff metric γ is more efficient than uniform FDM on clock metric in terms of convergence rate and computational time for fractional relaxation equation on the condition of similar accuracy.

**Table1.** The comparison between the implicit S-FDM and implicit uniform FDM on(17) when fractional derivative orders are α=0.5 and initial value A=10, relaxation parameter B=1, time terminal T=20 and space collocation nodes 100. CPU time is the computation cost on mentioned computer condition.

Maximum relative error(MRE) $= \max(abs(u(:,t_N) - \bar{u}(:,t_N))/u(:,t_N))$

Convergence rate $= \log_2\left(\dfrac{MRE(N/2,M)}{MRE(N,M)}\right)$

| Nodes of scales-FDM | MRE | Convergence rate | CPU time(s) |
|---|---|---|---|
| 25 | 0.0024 | | 0.320132 |
| 50 | 8.4268e-04 | 1.51 | 0.632755 |
| **100** | **2.9577e-04** | **1.5105** | **1.247940** |
| 200 | 1.0416e-04 | 1.5057 | 2.452703 |
| Nodes of uniform FDM | MRE | Convergence rate | CPU time(s) |
| 25 | 0.0103 | | 0.441325 |
| 50 | 0.0050 | 1.0426 | 0.869365 |
| 100 | 0.0025 | 1.0 | 1.726982 |
| 200 | 0.0012 | 1.0589 | 3.444368 |
| **400** | **6.0580e-04** | **0.9861** | **6.849005** |

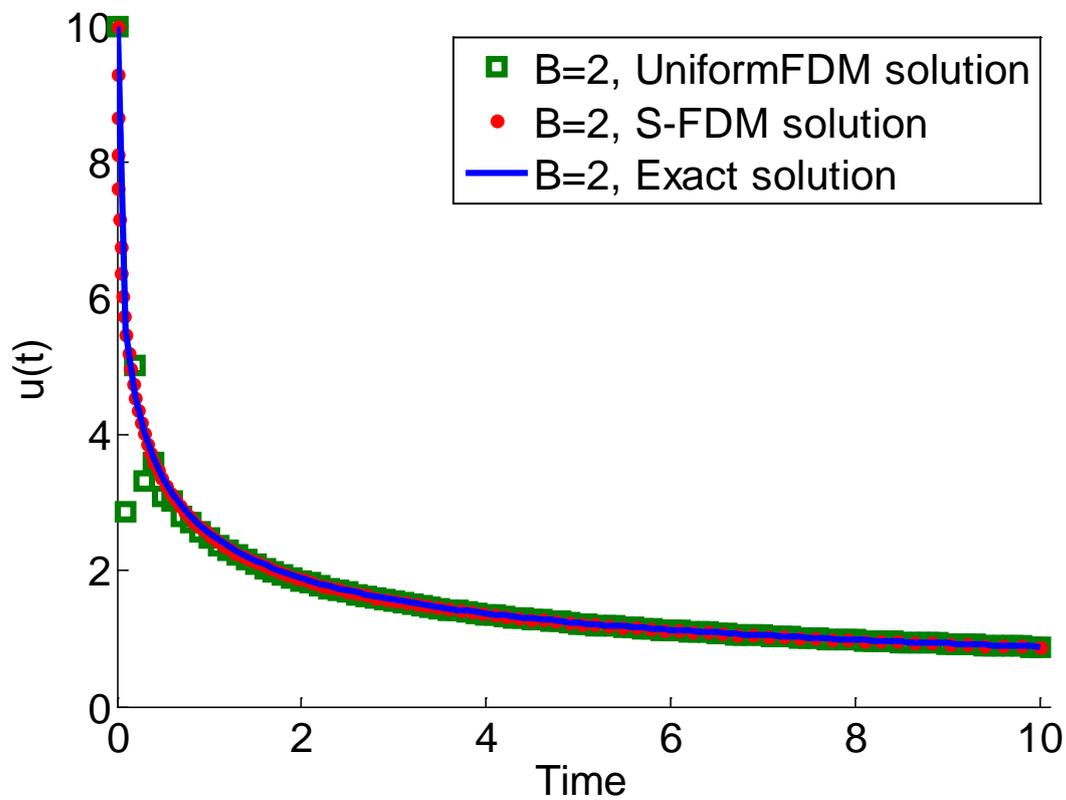

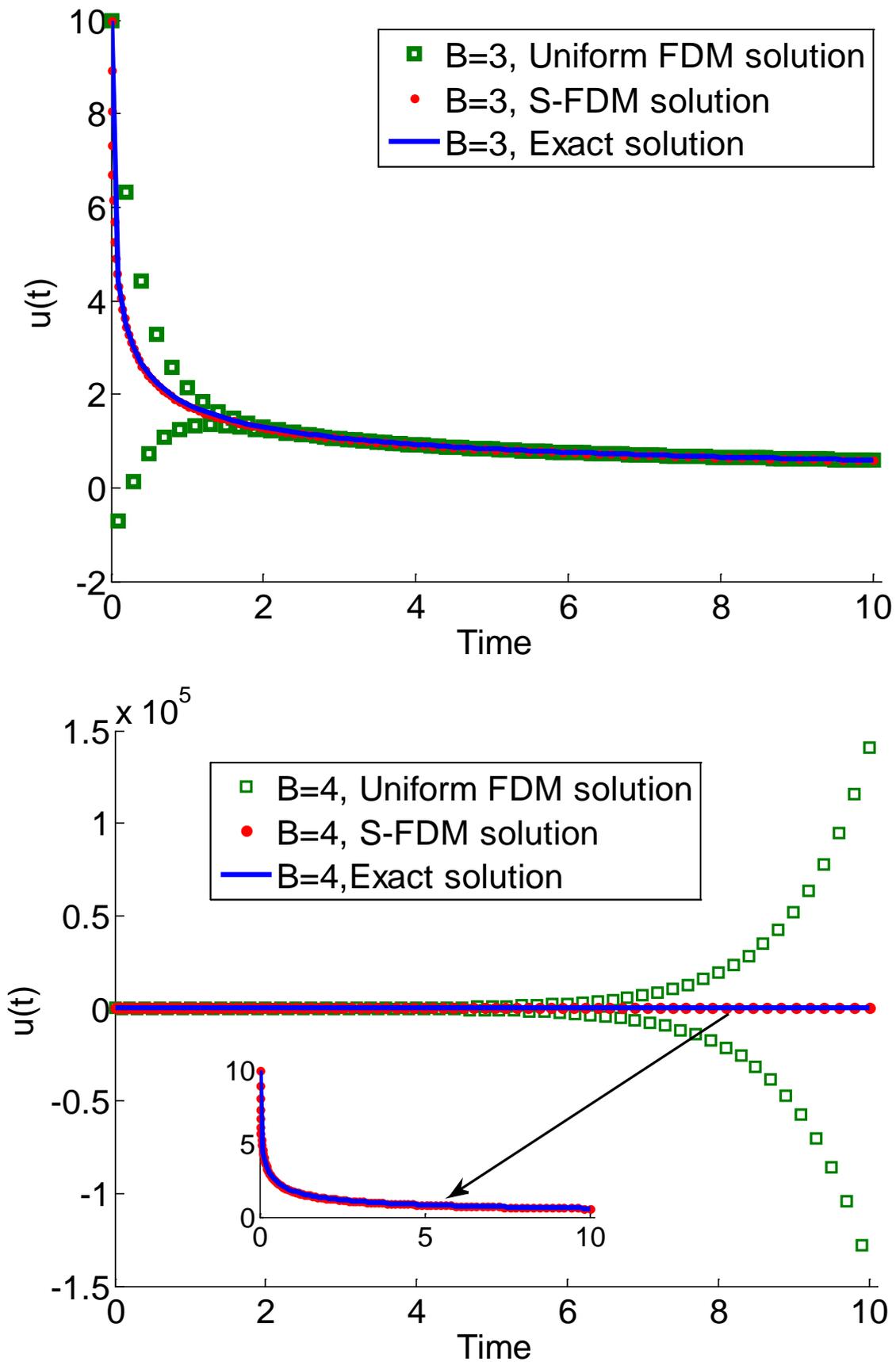

**Fig. 1.** Numerical results with γ=0.5 by S-FDM and uniform FDM.

**Example2.** Consider the following fractional diffusion equation

$$\begin{cases} \dfrac{\partial^\alpha u(x,t)}{\partial t^\alpha} = k \dfrac{\partial^2 u(x,t)}{\partial x^2}, \ x \in (0, L), \ t > 0 \\ u(x,0) = \sin\left(\dfrac{\pi x}{L}\right), \ x \in [0, L] \\ u(0,t) = 0, \ u(L,t) = 0 \end{cases} \quad (17)$$

If diffusion coefficient selects $k = L^2/\pi^2$, the exact solution is $u_{exact} = \sin(\pi x/L)E_\gamma(-t^\alpha)$.

On Table 2, a variety of calculation characteristics (accuracy, rate of convergence and computation cost) are compared between the two methods on the same condition. It reveals that the S-FDM for diffusion problems could nearly arrange one quarter collocations and increase the computation efficiency 7 times while the uniform FDM has 400 collocations and similar accuracy, separately. It also shows the proposed method appears a super-linear convergence rate (the rates are near 1.4(2-α)) while the uniform FDM just has a linear convergence. Therefore, the S-FDM reveals higher computation efficiency and accuracy. Meanwhile, for the sake of further illustration of computation property, we make Table 3 in which total time step of uniform FDM is equal to 0.1. It clearly shows the discrepancy of computation cost between S-FDM and uniform FDM while the numerical results reach the similar max relative error. It shows that S-FDM could use few nodes and improve computation efficient for long time range apparently, in especial S-FDM just expend ~1/400CPU time of uniform mesh method with better relative error for time T=500, when the fractional order α is equal to 0.4.

On Fig.2, the diffusion equation restarts with the natural (Neumann) boundary condition (gradient is equal to 0) and point initial condition ($u(L/2,0) = 2, u(\text{else},0) = 0$) on equation(17) when the space end L=10, time terminal T=10. From Fig.2(a), it indicates that different temporal collocations make completely similar solutions with S-FDM. In contrast, the curves on Fig.2(b) reveals distinctly diacritical results of varies of collocations, especially in the initial time domain. Comparing Fig.2(a) with Fig.2(b), the S-FDM in power-law time behaves better for diffusion problem with point source that we could choose few collocations to improve the efficiency with precise numerical results. Furthermore, the proposed numerical scheme could make numerical results more accurate and stable for diffusion equation.

**Table2.** The comparison between the implicit S-FDM and implicit uniform mesh method with when fractional derivative orders are α=0.6 and L=10, space collocation nodes 100, time terminal T=20. CPU time is the computation cost on mentioned computer condition. Here, MRE and Convergence rate have the same definition as the Table 1.

| Nodes of S-FDE | MRE | Convergence rate | CPU time(s) |
|---|---|---|---|
| 25 | 0.0050 | | 0.539049 |
| 50 | 0.0020 | 1.5850 | 1.272000 |
| **100** | **8.0701e-04** | **1.3093** | **3.799202** |
| 200 | 3.6461e-04 | 1.1462 | 10.141070 |
| Nodes of uniform FDM | MRE | Convergence rate | CPU time(s) |
| 25 | 0.0138 | | 0.631350 |
| 50 | 0.0067 | 1.0424 | 1.417968 |
| 100 | 0.0033 | 1.0217 | 3.448372 |
| 200 | 0.0017 | 0.9569 | 9.425320 |
| **400** | **8.7993e-04** | **0.9501** | **29.014906** |

**Table 3.** The comparison of computation cost between the implicit S-FDM and implicit uniform mesh method with when fractional derivative orders are α=0.4, 0.6, 0.8 seperately and L=10，space collocation 100，time step size of uniform mesh is 0.1. CPU time is the computation cost on mentioned computer condition. Here，MRE is same as the Table 1.

| T(α=0.4) | Nodes of S-FDM | MRE | CPU time | Nodes of uniform FDM | MRE | CPU time |
|---|---|---|---|---|---|---|
| 1 | 6 | 0.0128 | 0.132261 | 10 | 0.0134 | 0.253232 |
| 10 | 20 | 0.0019 | 0.464425 | 100 | 0.0018 | 3.411497 |
| 50 | 45 | 0.00043568 | 1.267172 | 500 | 0.00044247 | 41.671458 |
| 100 | 60 | 0.00026592 | 1.860014 | 1000 | 0.00026494 | 143.707843 |
| 200 | 80 | 0.00017538 | 2.856601 | 2000 | 0.00017397 | 529.405953 |
| 500 | 150 | 0.00010578 | 7.635067 | 5000 | 0.00011852 | 3075.40376 |
| T(α=0.6) | Nodes of S-FDM | MRE | CPU time | Nodes of uniform FDM | MRE | CPU time |
| 1 | 8 | 0.0231 | 0.179223 | 10 | 0.0248 | 0.234092 |
| 10 | 40 | 0.0031 | 1.088697 | 100 | 0.0034 | 3.341204 |
| 50 | 100 | 0.00064047 | 4.005253 | 500 | 0.00071016 | 40.793088 |
| 100 | 200 | 0.00026286 | 12.613004 | 1000 | 0.00039231 | 141.362693 |
| 200 | 350 | 0.00015362 | 34.449055 | 2000 | 0.00023614 | 523.511651 |
| 500 | 700 | 0.0001051 | 124.781522 | 5000 | 0.00014357 | 3094.11699 |
| T(α=0.8) | Nodes of S-FDM | MRE | CPU time | Nodes of uniform FDM | MRE | CPU time |

| 1 | 10 | 0.0369 | 0.238812 | 10 | 0.0404 | 0.235752 |
| 10 | 70 | 0.0064 | 2.355889 | 100 | 0.0064 | 3.372816 |
| 50 | 250 | 0.00088249 | 18.758857 | 500 | 0.00096826 | 41.596482 |
| 100 | 400 | 0.00088249 | 43.725121 | 1000 | 0.00049561 | 143.212046 |
| 200 | 700 | 0.00027866 | 124.939695 | 2000 | 0.00028934 | 528.658683 |
| 500 | 1500 | 0.00015656 | 550.187193 | 5000 | 0.00016396 | 3096.849746 |

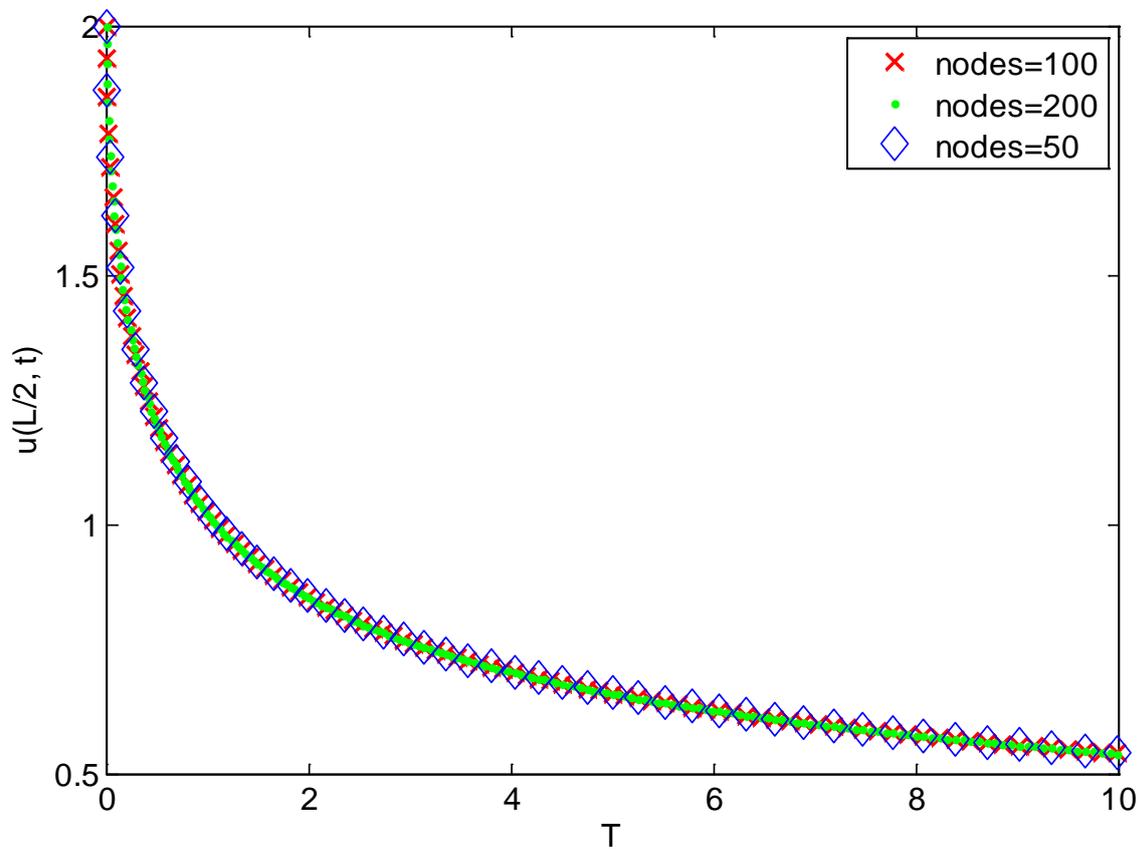

(a)

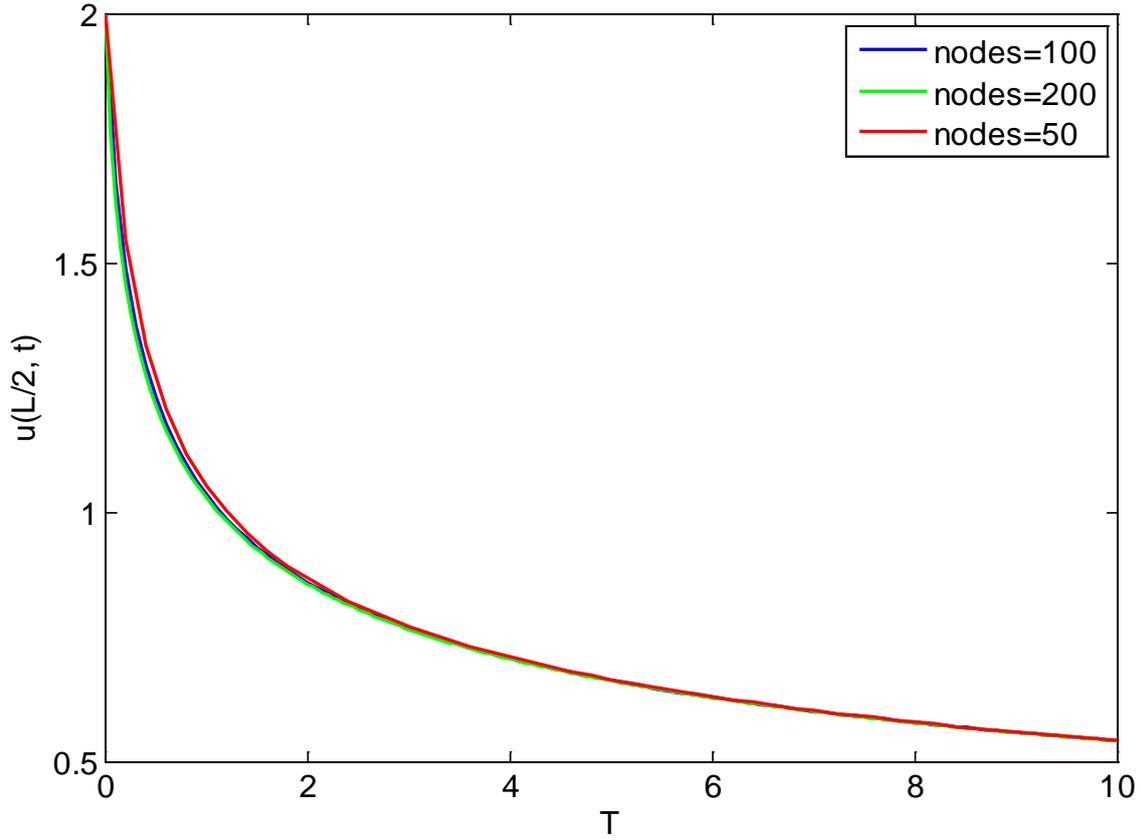

(b)

**Fig. 2.** Numerical result of S-FDM is shown in figure (a), while figure (b) adopts the uniform mesh method. Both of them take implicit difference scheme on condition of diffusion coefficient D=0.005 and fractional order α=0.6 with collocation nodes 50, 100, 200 at time.

**Example3.** Consider the point source type fractional diffusion equation of the form

$$\begin{cases} \dfrac{\partial^\alpha u(x,t)}{\partial t^\alpha} = k\dfrac{\partial^2 u(x,t)}{\partial x^2}, \; x \in (0,L), \; t>0 \\ u(x,0) = \begin{cases} 2 & x=L/2 \\ 0 & x \neq L/2 \end{cases} \\ u(0,t)=0, \; \dfrac{\partial u(0,t)}{\partial u}=0 \end{cases} \quad (18)$$

where the diffusion coefficient is 0.005. Here, we regard $\bar{u}(x,t_{N=5000})$ as the approximation exact solution while the time collocations of numerical solution $\bar{u}(x,t_{N\leq 100})$ is much less than collocation nodes of approximation exact solution $\overline{u}(x, t_{N=5000})$.

On Table 4, numerical solutions for the fractional diffusion equation(18) are presented. Here, three collocation methods are carried out in which the condition of non-uniform meshes satisfies the relationship [14, 28] $\tau_n = (N+1-n)\mu$ where $\mu = \dfrac{2T}{N(N+1)}$. The accuracy and convergence rate are compared with a variety of fractional orders and collocation nodes which indicated that the S-FDM could improve the numerical accuracy apparently and shows (2-α) rate of convergence. Meanwhile, it is observed that the traditional uniform FDM and non-uniform FDM mentioned gain worse accuracy and just 1-st convergence rate for the source type diffusion equation. Furthermore, through comparison on the accuracy between the uniform mesh method and non-uniform method mentioned, the conclusion that the non-uniform method behaves worse for the distinct relaxation diffusion equation could be drawn.

**Table4.** The comparison between the implicit S-FDM and implicit uniform FDM when fractional derivative orders α are equal to 0.8,0.6,0.4, separately on (18) with the condition L=10，space collocation 100，time terminal T=10. CPU time is the computation cost on mentioned computer condition.

Maximum absolute error $1(\text{MAE1}) = \max(abs(\overline{u}(:,t_{N=5000}) - \overline{u}(:,t_{N\ll 5000})))$

$\text{Rate1} = \log_2\left(\dfrac{\text{MAE1}(N/2, M)}{\text{MAE1}(N, M)}\right)$

| S-FDM | MAE1(α=0.8) | Rate1 | MAE1 (α=0.6) | Rate1 | MAE1 (α=0.4) | Rate1 |
|---|---|---|---|---|---|---|
| Nodes=20 | 0.0041 | | 0.0022 | | 7.7209e-04 | |
| Nodes=40 | 0.0018 | 1.1876 | 8.3302e-04 | 1.4011 | 2.4282e-04 | 1.6689 |
| Nodes=80 | 8.3141e-04 | 1.1144 | 3.1415e-04 | 1.4069 | 7.0857e-05 | 1.7769 |
| Uniform FDM | MAE1(α=0.8) | Rate1 | MAE1 (α=0.6) | Rate1 | MAE1 (α=0.4) | Rate1 |
| Nodes=20 | 0.0049 | | 0.0044 | | 0.0034 | |
| Nodes=40 | 0.0024 | 1.0297 | 0.0021 | 1.0671 | 0.0017 | 1.0 |
| Nodes=80 | 0.0011 | 1.1255 | 0.0010 | 1.0704 | 8.1807e-04 | 1.0552 |
| Non-uniform FDM | MAE1(α=0.8) | Rate1 | MAE1 (α=0.6) | Rate1 | MAE1 (α=0.4) | Rate1 |
| Nodes=20 | 0.0079 | | 0.0079 | | 0.0064 | 1 |
| Nodes=40 | 0.0039 | 1.0184 | 0.0039 | 1.0184 | 0.0032 | 1 |
| Nodes=80 | 0.0020 | 0.9635 | 0.0020 | 0.9635 | 0.0016 | 1 |

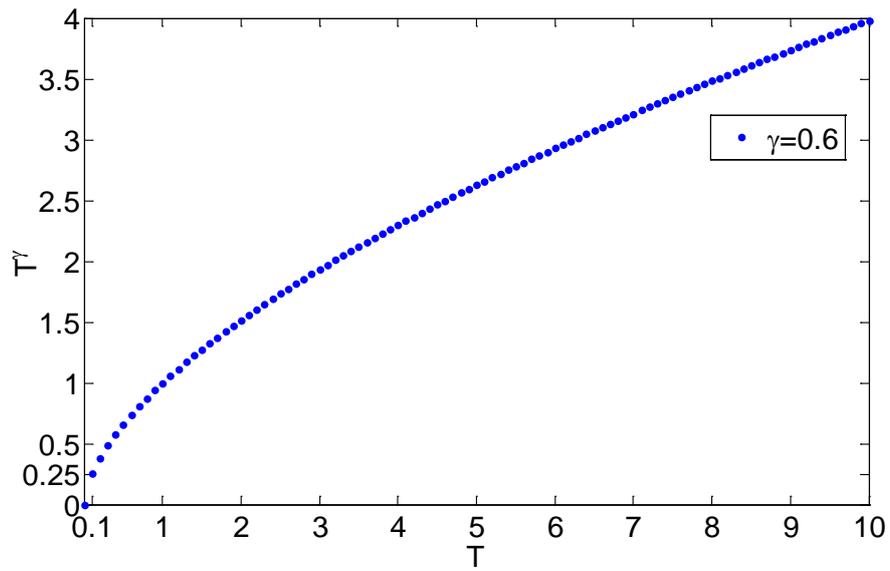

(a)

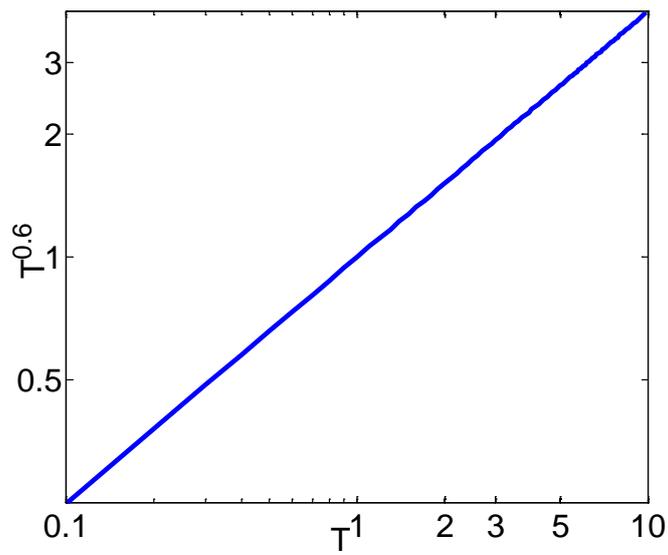

(b)

**Fig. 3.** The point positions on power-law time $T^{0.6}$ (Hausdroff metric **γ**=0.6) on contrast parallel positions on clock time when collocation nodes are 100.

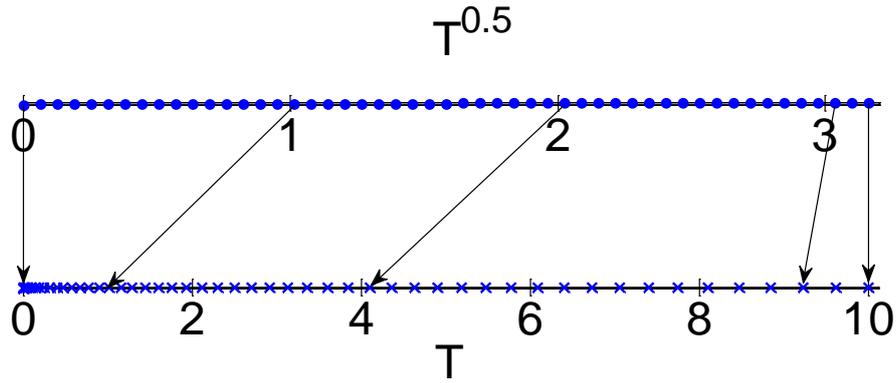

**Fig. 4.** The point positions on power-law time (Hausdroff metric **γ**=0.5) with the equidistant step and the parallel positions on clock time when collocation nodes are 50.

### 4. Discussion

The proposed method establishes the relationship fractional derivative order and strategy of non-uniform collocation mode by scale transfer based on Hausdroff metric. It means the proposed method behaves uniform mesh on power-law time scale, while it shows non-uniform mesh on clock time scale. It is clear shown that S-FDM collocates equidistant step-size nodes on power-law time $T^\alpha$ in Fig. 4. This kind of collocation mode leads dense nodes at early time and sparse nodes at later time on clock time scale. On contrast, the step-sizes are larger than others at the beginning of the nodes on power-law time scale $T^\gamma$ when the mesh is equidistant step-sizes on clock time scale T in Fig.3(a). For example, time is equal to 0.1 in clock time while the time value becomes 0.25 in power-law time. This situation leads to numerical results unstable and loose in numerical simulation which is shown on Figures (1-2) from numerical experiments. Particularly, it may cause a lager calculation error in the initial time.

To verify our statement, here fractional derivative relaxation equation is considered. In these examples, the S-FDM keeps the (2-α)-th order convergence, while the uniform FDM and non-uniform FDM [14, 28]lose their convergence rate and accuracy from the numerical results while only keep 1-th order convergence rate. In addition, it is easy to conclude that the proposed method can improve the stability of numerical results apparently. For example, Fig 2 shows the S-FDM allows the wide range variation of relaxation coefficient (*B*).

in Eq. (1), compared with FDM. Moreover, the S-FDM with easy-to-implement collocation strategy (Appendix) can achieve higher accuracy compared with uniform FDM, which benefits the programming and computation cost. The high computation efficiency of S-FDM has been clearly shown in Table 3, especially for low order time FREs and long time range computation. For example, S-FDM just need ~1/400CPU time of uniform mesh method with better relative error at time $T$=500 when the fractional order α is equal to 0.4. If we use identical collocation nodes (i.e., $N$=150), the accuracy of S-FDM is much better than FDM, for example, the accuracy of S-FDM (MRE is 0.00010578) is up to one order of magnitude higher than FDM (MRE is 0.0014) for $N$=150.

## 5. Conclusion

The S-FDM can be regarded as a uniform mesh method on the view of Hausdroff metric whilst non-uniform mesh method from the viewpoint of Euclidean metric. In this paper, we proved that the proposed method is unconditional stable with (2-α)-th convergence rate. Meanwhile, numerical experiments clearly demonstrate the advantage of the S-FDM in comparison with uniform mesh FDM and non-uniform FDM [14, 28] in terms of the accuracy, rate of convergence and computation cost. It can be concluded that, the S-FDM offers an easy-to implement and roust non-uniform collocation mode based on scale metric transfer approach, to efficiently solve FREs. Though, we only considered the time FREs, the S-FDM can be extended to various kinds of time and space fractional differential equations by applying the relationship between fractional derivative order and strategy of non-uniform collocation.

**Acknowledgment**

The work was supported by the National Natural Science Foundation of China (Grant Nos. 41330632, 11572112, 41628202), the Postgraduate Research & Practice Innovation Program of Jiangsu Province (Grant Nos.KYCX17_0488, KYCX17_0490) and the Fundamental Research Funds for the Central Universities (2017B710X14, 2017B709X14).

**Appendix**

**1 Formulation in numerical schemes of relaxation equation and diffusion equation**

Anomalous relaxation and diffusion is commonly observed in field experiment, especially in unsaturated flow and transportation on underground water. To describe the memory effect of anomalous diffusion phenomena, the Caputo time fractional relaxation equation is given by

$$\begin{cases} {}_0^C D_t^\gamma u(t) + Bu(t) = f(t) & 0 \leq n-1 < \gamma \leq n \leq 2 \\ u^{(k)}(0) = u_0^{(k)} & (k = 0,1,\ldots,n-1) \end{cases} \quad (19)$$

where $0 < \gamma \leq 1$ represents fractional relaxation equation.

At the same, the time fractional diffusion equation is

$$\begin{cases} \dfrac{\partial^\alpha u(x,t)}{\partial t^\alpha} = D(x,t)\dfrac{\partial^2 u(x,t)}{\partial x^2} + f(x,t) \\ u(x,t) = h(x,t), x \in \partial\Omega, \ t > 0, \\ u(x,0) = w(x,t), \ t = 0. \end{cases} \quad (20)$$

where $\alpha$ is the fractional order and $\alpha \in (0,1]$ and $x \in [0,L]$, $t > 0$, D(x,t) is the diffusion parameter with dimension of $[L^2 T^{-\alpha}]$ and $f(x,t)$ represents source term.

As is known to us, the initial problems(19) could convert to the Volterra integral equation.

$$u(t) = \sum_{i=0}^{n-1} u_0^{(j)} \frac{t^j}{j!} + \frac{1}{\Gamma(\gamma)} \int_0^t (t-\tau)^{\gamma-1} (f(\tau) - Bu(\tau))d\tau \quad (21)$$

Here, let $\Delta_t u(x,t_k) = u(x,t_{k+1}) - u(x,t_k)$, an $db_{k,j} = k^{1/\alpha} - j^{1/\alpha}$,

we could give a numerical scheme of rectangular quadrature formula on the power-law time

$$\int_0^{t_{k+1}} (t_{k+1}-\tau)^{\gamma-1} (f(\tau) - Bu(\tau))d\tau \approx \sum_{j=0}^{k} m_{j,j+1}(f(t_j) - Bu_j) \quad (22)$$

where $m_{j,k} = \dfrac{\Delta t_\gamma}{\gamma}(b_{n+1,j}^\gamma - b_{n+1,k+1}^\gamma)$

and we could have a lower triangular matrix $\mathbf{M}_{(n+1)*(n+2)}$ as follow

$$\mathbf{M} = \begin{bmatrix} b_{1,0} & b_{1,1} & 0 & \cdots & 0 & 0 \\ b_{2,0} & b_{2,1} & b_{2,2} & \cdots & 0 & 0 \\ \vdots & & \ddots & \ddots & & \vdots \\ b_{n,0} & b_{n,1} & b_{n,2} & \cdots & b_{n,n} & 0 \\ b_{n+1,0} & b_{n+1,1} & b_{n+1,2} & \cdots & b_{n+1,n} & b_{n+1,n+1} \end{bmatrix}_{(n+1)\times(n+2)}$$

Here, we also could solve the relaxation equation through the implicit scheme after scale-dependent difference of the Caputo time fractional derivative

$$\frac{\partial^{\gamma} u(t_{n+1})}{\partial t^{\gamma}} \approx \frac{1}{\Gamma(1-\gamma)} \sum_{k=0}^{n} \frac{\Delta_t u(t_k)}{\tau_{k+1}} \int_{t_k}^{t_{k+1}} \frac{d\xi}{(t_{n+1}-\xi)^{\gamma}}$$
$$= \frac{-\Delta t_{\gamma}^{-1}}{\Gamma(2-\alpha)} \sum_{k=0}^{n} \frac{\Delta_t u(t_k)}{b_{k+1,k}} [(b_{n+1,k+1})^{1-\gamma} - (b_{n+1,k})^{1-\gamma}] \quad (23)$$

where $\tau_k = t_k - t_{k-1}$. The $u(t)$ is the solution of the equation(1).

Combining the implicit difference scheme and equations(9), the time fractional relaxation equation with implicit scheme can possess the numerical results with S-FDM as follow

n=0
$$\left(\frac{\Delta t_{\alpha}^{-1}}{\Gamma(2-\alpha)} + B\right) u(t_1) = \frac{\Delta t_{\alpha}^{-1}}{\Gamma(2-\alpha)} u(t_0) + f(t_1)$$

n=1
$$\left(\frac{\Delta t_{\alpha}^{-1}}{\Gamma(2-\alpha)} + B\right) u(t_2) = f(t_2) + \frac{\Delta t_{\alpha}^{-1}}{\Gamma(2-\alpha)} \left(\left(\frac{b_{2,1}^{1-\alpha}}{b_{2,1}} + b_{2,1}^{1-\alpha} - b_{2,0}^{1-\alpha}\right) u(t_1) - \left(b_{2,1}^{1-\alpha} - b_{2,0}^{1-\alpha}\right) u(t_0)\right)$$

n=2
$$\left(\frac{\Delta t_{\alpha}^{-1}}{\Gamma(2-\alpha)} \frac{b_{n+1,n}^{1-\alpha}}{b_{n+1,n}} + B\right) u(t_{n+1}) =$$

$$\frac{\Delta t_{\alpha}^{-1}}{\Gamma(2-\alpha)} \left( \begin{array}{c} \left(\frac{b_{n+1,n}^{1-\alpha}}{b_{n+1,n}} + \frac{b_{n+1,n}^{1-\alpha} - b_{n+1,n-1}^{1-\alpha}}{b_{n,n-1}}\right) u(t_n) - \left(b_{n+1,1}^{1-\alpha} - b_{n+1,0}^{1-\alpha}\right) u(t_0) \\ + \sum_{k=1}^{n-1} u(t_k) \left(-\frac{b_{n+1,k+1}^{1-\alpha} - b_{n+1,k}^{1-\alpha}}{b_{k+1,k}} + \frac{b_{n+1,k}^{1-\alpha} - b_{n+1,k-1}^{1-\alpha}}{b_{k,k-1}}\right) \end{array} \right) + f(t_{n+1})$$

where $x_i$ represents the space position and $t_i$ stands for time position.

Compared with the numerical scheme on clock time, the S-FDM appears to be complex and high computation efficient due to its complexity formulas which is relative with part $b_{k,j}$. However, the complexity formulas show interesting generalizations to accelerate the computation efficiency by reading a lower triangular matrix $\mathbf{M}_{(n+1)*(n+1)}$ and matrix $\mathbf{M_1}=\mathbf{M}.\wedge(1-\alpha)$ and make the results more precision.

$$\mathbf{M} = \begin{bmatrix} b_{1,0} & 0 & \cdots & 0 & 0 \\ b_{2,0} & b_{2,1} & \cdots & 0 & 0 \\ \vdots & & \ddots & & \vdots \\ b_{n,0} & b_{n,1} & \cdots & b_{n,n-1} & 0 \\ b_{n+1,0} & b_{n+1,1} & \cdots & b_{n+1,n-1} & b_{n+1,n} \end{bmatrix}_{(n+1)\times(n+1)}$$

The time fractional diffusion equation is

n=0

$$\left(\frac{\Delta t_\alpha^{-1}}{\Gamma(2-\alpha)}+\frac{2}{\Delta x^2}\right)u(x_i,t_1)-D\frac{u(x_{i+1},t_1)+u(x_{i-1},t_1)}{\Delta x^2}=\frac{\Delta t_\alpha^{-1}}{\Gamma(2-\alpha)}u(x_i,t_0)+f(x_i,t_1)$$

n=1

$$\left(\frac{\Delta t_\alpha^{-1}}{\Gamma(2-\alpha)}+\frac{2}{\Delta x^2}\right)u(x_i,t_2)-D\frac{u(x_{i+1},t_2)+u(x_{i-1},t_2)}{\Delta x^2}=f(x_i,t_2)+$$

$$\frac{\Delta t_\alpha^{-1}}{\Gamma(2-\alpha)}\left(\left(\frac{b_{2,1}^{1-\alpha}}{b_{2,1}}+b_{2,1}^{1-\alpha}-b_{2,0}^{1-\alpha}\right)u(x_i,t_1)-\left(b_{2,1}^{1-\alpha}-b_{2,0}^{1-\alpha}\right)u(x_i,t_0)\right)$$

n=2

$$\left(\frac{\Delta t_\alpha^{-1}}{\Gamma(2-\alpha)}\frac{b_{n+1,n}^{1-\alpha}}{b_{n+1,n}}+\frac{2}{\Delta x^2}\right)u(x_i,t_{n+1})-D\frac{u(x_{i+1},t_{n+1})+u(x_{i-1},t_{n+1})}{\Delta x^2}=$$

$$\frac{\Delta t_\alpha^{-1}}{\Gamma(2-\alpha)}\left(\begin{array}{l}\left(\dfrac{b_{n+1,n}^{1-\alpha}}{b_{n+1,n}}+\dfrac{b_{n+1,n}^{1-\alpha}-b_{n+1,n-1}^{1-\alpha}}{b_{n,n-1}}\right)u(x_i,t_n)-\left(b_{n+1,1}^{1-\alpha}-b_{n+1,0}^{1-\alpha}\right)u(x_i,t_0)\\ +\displaystyle\sum_{k=1}^{n-1}u(x_i,t_k)\left(-\dfrac{b_{n+1,k+1}^{1-\alpha}-b_{n+1,k}^{1-\alpha}}{b_{k+1,k}}+\dfrac{b_{n+1,k}^{1-\alpha}-b_{n+1,k-1}^{1-\alpha}}{b_{k,k-1}}\right)\end{array}\right)+f(x_i,t_{n+1})$$

where $x_i$ represents the space position and $t_i$ stands for time position.